\documentclass{acm_proc_article-sp}

\newtheorem{theorem}{Theorem}
\newtheorem{prop}{Proposition}
\newtheorem{conj}{Conjecture}
\newdef{example}{Example}
\newdef{question}{Question}

\newcommand{\rr}{{\mathbb R}}

\newcommand{\ba}{{\mathbf a}}

\newcommand{\sA}{{\mathcal A}}

\begin{document}

\title{Toric Ideals of Homogeneous Phylogenetic Models}

\numberofauthors{1}
\author{
  \alignauthor Nicholas Eriksson\\
  \affaddr{Department of Mathematics}\\
  \affaddr{University of California, Berkeley}\\
  \affaddr{Berkeley, CA 94720-3840}
  \email{eriksson@math.berkeley.edu}
}

\date{January 7, 2004}
\maketitle

\begin{abstract}
We consider the model of 
phylogenetic trees in which every node of the tree is an observed,
binary random variable and the transition probabilities are given by
the same matrix on each edge of the tree.
The ideal of invariants of this model is a toric ideal in 
$\mathbb C[p_{i_1 \dots i_n}]$.
We are able to compute the Gr\"obner basis and minimal generating set
for this ideal for trees with up to 11 nodes. 
These are the first non-trivial Gr\"obner bases calculations in
$2^{11} =  2048$ indeterminates.  We conjecture that there is a quadratic
Gr\"obner basis for binary trees, but that generators of degree $n$ 
are required for some trees with $n$ nodes.
The polytopes associated with these toric ideals display interesting
finiteness properties.  We describe the polytope for an infinite
family of binary
trees and conjecture (based on extensive computations) that there is a
universal bound on the number of vertices of the polytope of a binary
tree.
\end{abstract}

\section{Introduction}\label{sec:intro}

%
%
%
%
%
%
%




A phylogenetic tree is a rooted tree $T$ on $n$ nodes with a
$\kappa$-ary random variable $X_i$ associated to every node.  Write
$\rho(v)$ for the parent of node $v$.  Then the transition
probabilities between $\rho(v)$ and $v$ are given by a $\kappa$ by
$\kappa$ matrix $A^{(v)}$ for every non-root node of $T$.


In an application, 
$\kappa$ might encode the four nucleic acids that make up DNA, 
the two families of nucleic acids, or the twenty amino acids.
The transition matrices are generally picked from some specific family
such as the Jukes-Cantor \cite{fu}, Kimura \cite{steel}, or general Markov models \cite{allman}.

In this paper we consider the homogeneous Markov model where all $A^{(v)}$ are equal, all
nodes are binary ($\kappa = 2$) and observable, and the root has uniform distribution.  
We write
$A^{(v)} = A = \begin{pmatrix} a_{00} & a_{01} \\ a_{10} & a_{11}\end{pmatrix}$.
The probability of observing $i$ at a node $v$ is computed from
the parent of $v$ by
\[
P(X_v = i) = a_{0i}P(X_{\rho(v)} = 0) +  a_{1i}P(X_{\rho(v)} = 1).
\]
We are interested in the algebraic relations satisfied by the joint
distribution
\[
p_{i_1i_2\dots i_n} := P(X_1 = i_1, \dots, X_n = i_n).
\]
Writing the joint distribution in terms of the model parameters
$a_{00}, a_{01}, a_{10}, a_{11},$ we have
\begin{equation}\label{eq:toricdef}
p_{i_1i_2\dots i_n} = \prod_{j=2}^{n} a_{i_{\rho(j)} i_j}
\end{equation}
where the nodes of the tree are labeled $1$ to $n$ starting with the
root.  That is, the probability of observing a certain labeling
of the tree is the product of the $a_{ij}$ that correspond to the
transitions on all edges of the tree.  The indeterminates $a_{ij}$
parameterize a toric variety of dimension 4 in $\mathbb R^{2^n}$.
We let $I_T$ be the corresponding toric ideal, called 
the ideal of phylogenetic invariants.
In the notation of \cite{gbcp}, the toric ideal $I_T$ is specified by the 4
by $2^n$ configuration $\sA_T$, where column $(i_1, \dots, i_n)$ consists of
the exponent vector of the $a_{ij}$ in (\ref{eq:toricdef}).  We order
the rows $(a_{00}, a_{01}, a_{10}, a_{11})$.  Let $P_T$ be the convex hull
of the columns of $\sA_T$.

%
%
%
%
%
%
%

We are interested in two questions from \cite{bio}.
First, which relations on the joint probabilities $p_{i_1\dots i_n}$ does
the model imply?  This problem is solved by giving generators of
the ideal of invariants $I_T$.

In Section~2, we study the generators of this ideal.
Our main
accomplishment is the computation of Gr\"obner and Markov bases for
trees with 11 nodes.  These are computations in 2048 indeterminants,
which we believe to be the largest number of indeterminants ever in a
Gr\"obner basis calculation.  We also calculate generating sets for
all trees on at most 9 nodes.  Based on this evidence, we conjecture
that if $T$ is binary, then the ideal $I_T$ has a quadratic generating
set, and furthermore, that relations of degree $n$ are necessary to
generate $I_T$ for certain trees  with $n$ nodes.


Our second goal is to determine, given a labeling of the tree $T$, if
we can identify parameters $a_{ij}$ such that the labeling is the most
likely among all labelings?  This problem is solved by computing the
normal fan of the toric variety in the sense of \cite{fulton}.

In Section~3, we study this normal fan and the polytope $P_T$.  Our
main result, Theorem~\ref{thm:1}, is an explicit description of 
the polytope $P_T$ for an infinite family of binary trees.  
For this family, $P_T$ always has 8 vertices and 6 facets which we characterize.
We also present extensive calculations of $P_T$ for various trees
and conjecture that there is a bound on the number of vertices of $P_T$ 
as $T$ ranges over all binary trees.



The invariants vanish for a given distribution
$(p_{i_1 \dots i_n})$ essentially when that distribution comes from our
model.  
Thus the knowledge of the generators of this ideal is potentially 
very useful for fitting biological sequence data to a phylogenetic 
tree, as first noted by Cavender and Felsenstein \cite{cavender}.
While there has been much progress towards finding the ideal of
invariants for other phylogenetic models (see \cite{allman}, \cite{fu},
\cite{steel}),
the homogeneous model is
particularly attractive because the low number of parameters makes it
possible to compute non-trivial examples.  
Hopefully we can use the
homogeneous model to approximate in some sense 
the general model, perhaps by subdividing edges of the tree.

\begin{example}\label{ex:1}
Let $T$ be a path with 3 nodes.  Then
\[
\sA_T = 
\begin{pmatrix}
2& 1& 0& 0& 1& 0& 0& 0\\
0& 1& 1& 1& 0& 1& 0& 0\\
0& 0& 1& 0& 1& 1& 1& 0\\
0& 0& 0& 1& 0& 0& 1& 2
\end{pmatrix},
\]
the polytope $P_T$ has 7 vertices and 6 facets,
and the toric ideal of the path of length 3 is generated by 6 binomials
\begin{equation*}
\begin{split}
I_T =  \langle &  x_{101}-x_{010},  x_{001}x_{100}-x_{000}x_{010},\\ 
& x_{011}x_{100}-x_{001}x_{110}, x_{011}x_{110}-x_{010}x_{111}, \\
& x_{001}^2x_{111}-x_{000}x_{011}^2,  x_{100}^2x_{111}-x_{000}x_{110}^2 \rangle.
\end{split}
\end{equation*}
\end{example}

\section{Toric Ideals}\label{sec:toric}


The toric ideals $I_T$ are homogeneous, since all monomials in
(\ref{eq:toricdef}) have the same degree $n-1$.  Thus they define
projective toric varieties $Y_T$.  Algebraic geometers usually require
a toric variety to be normal, but the reader should be warned that the
toric varieties discussed in this paper are generally not normal.

Recall that a projective toric variety given by a configuration $\sA = (\ba_1, \dots, \ba_k)$
is covered by the affine toric varieties given by $\sA -  \ba_i$.
An affine toric variety defined by a configuration $\sA$ is said to be smooth if 
the semigroup $\mathbb N \sA$ is isomorphic to $\mathbb N^r$ for some $r$ 
\cite[Lemma 2.2]{toric}.
\begin{prop}
The projective toric variety $Y_T$ of a binary tree $T$ is not smooth.
\end{prop}
\begin{proof}
Recall that the columns of the configuration $\sA_T$ are indexed by
$0/1$-labelings of the tree $T$.  Look at the affine chart $I_{\sA - \ba_{0\dots 0}}$,
where $\ba_{0\dots 0}$ corresponds to the all zero tree.
On this chart, write $\tilde{\ba}_{i} = \ba_{i} - \ba_{0\dots 0}$.
Let $10\dots 0$ be the tree with a 1 at the root and zeros everywhere else,
$0\dots 01$ be the tree with a 1 at a single leaf and zeros everywhere else,
and $0\dots 010\dots 0$ be the tree with a 
single 1 at the parent of a leaf and zeros elsewhere.
Then since $\ba_{0\dots 0} = (n-1,0,0,0)$, we have
\begin{gather*}
\tilde{\ba}_{10\dots 0} = (n-3,0,2,0) -\ba_{0\dots 0} = (-2,0,2,0)\\
\tilde{\ba}_{0\dots 01} = (n-2,1,0,0) -\ba_{0\dots 0}  = (-1,1,0,0)\\
\tilde{\ba}_{0\dots 010\dots 0} = (n-4,1,2,0) - \ba_{0\dots 0} = (-3,1,2,0),
\end{gather*}
and so we see that
\[
\tilde{\ba}_{10\dots 0} + \tilde{\ba}_{0\dots 01}  = \tilde{\ba}_{0\dots 010\dots 0}.
\]
Therefore, $\sA - \mathbf a_{0\dots 0}$ is not isomorphic to $\mathbf
N^r$ and the toric variety $Y_T$ is not smooth.
\end{proof}

We are primarily interested in the generators of the ideals $I_T$.
Knowledge of the generators would allow us to easily compute whether
given data came from the homogeneous Markov model from some specific
phylogenetic tree.

Using {\tt 4ti2} \cite{4ti2}, Gr\"obner and Markov
bases for the ideal $I_T$ were computed for all trees with at most 9
nodes as well as selected trees with 10 and 11 nodes.  This took about
6 weeks of computer time in total on a 2GHz computer.  The
computations in 2048 variables (trees with 11 nodes) each took as long as a
week and required over 2 GB of memory.

Details about the Markov bases for all binary trees with at most 11
nodes are shown in 
Table~1.
These computations lead
us to make the following conjectures.
\begin{table}
\begin{center}
\begin{tabular}{l|l|l|l}
tree & Degree & \#Minimal & Max degree \\
     & of $I_T$       & Generators & of generator\\\hline
\psfig{file=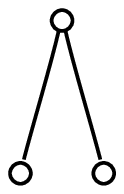, height=.13in}& 4 & 4 & 2 \\
\psfig{file=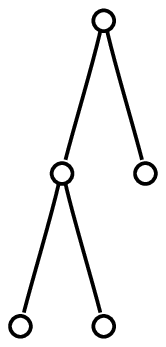, height=.24in}& 28 & 79 & 2 \\
\psfig{file=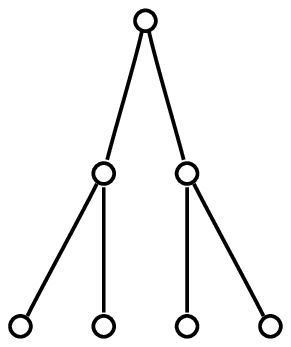, height=.24in}& 92 & 441 & 2 \\
\psfig{file=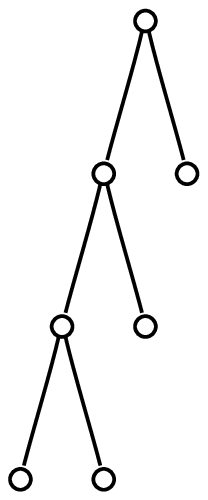, height=.35in}& 96 & 561 & 2 \\
\psfig{file=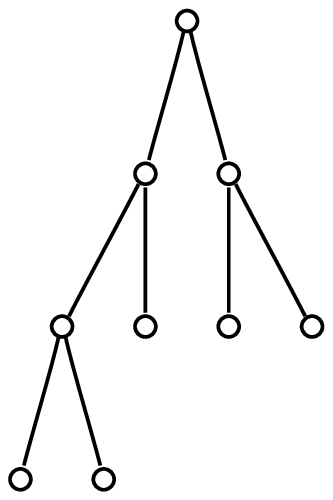, height=.35in}& 210 & 2141 & 2 \\
\psfig{file=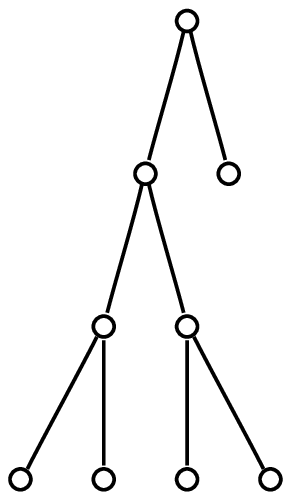, height=.35in} & 220 & 2068 & 2 \\
\psfig{file=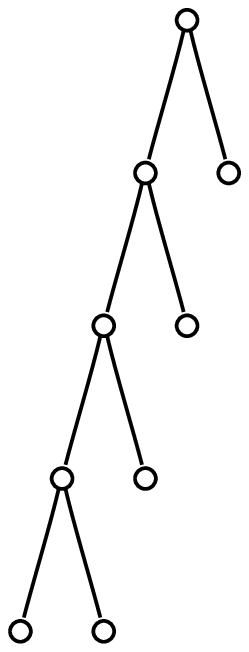, height=.42in} & 210 & 2266 & 2 \\
\psfig{file=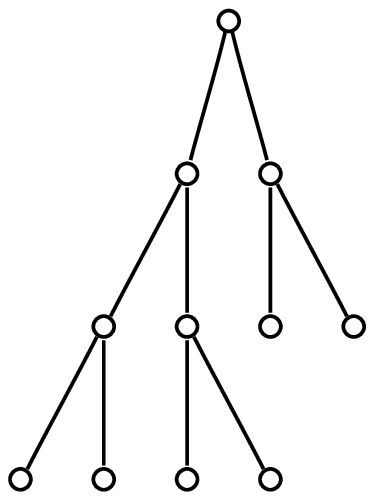, height=.35in} & 412 & 7121 & 2\\
\psfig{file=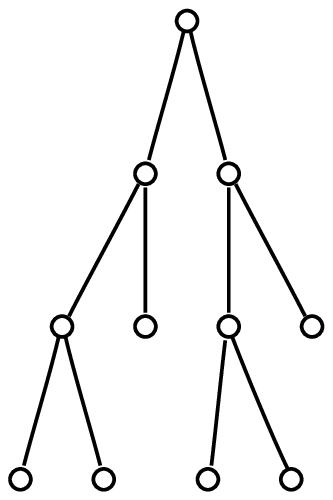, height=.35in} & 404 & 7131 & 2\\
\psfig{file=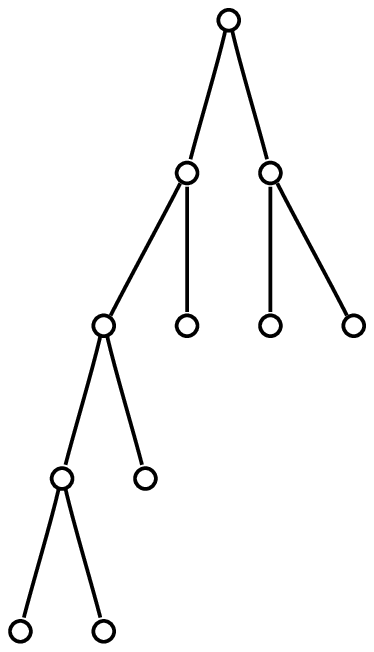, height=.42in}  & 400 & 7137 & 2\\
\psfig{file=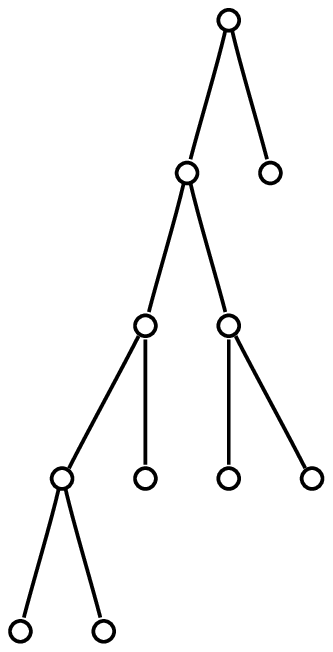, height=.42in}  & 412 & 7551 & 2\\
\psfig{file=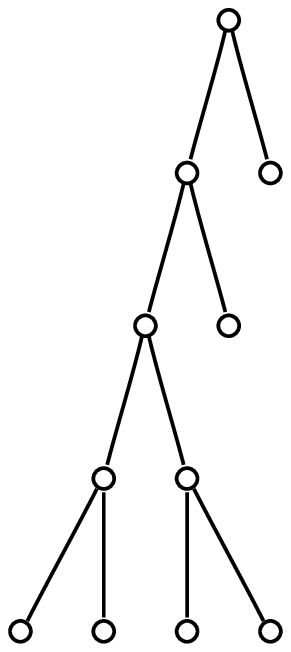, height=.42in}  & 412 & 7551 & 2\\
\psfig{file=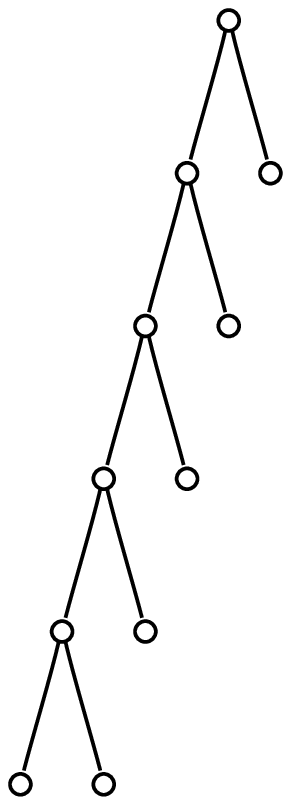, height=.55in}  & 404 & 7561 & 2\\
\end{tabular}
\end{center}
\caption{Degree of $I_T$, number of minimal generators, and maximum
degree of the generators}
\label{tab:binmarkov}
\end{table}
\begin{conj}
The toric ideal corresponding to a binary tree is generated in degree 2.
More generally, if every non-leaf node of the tree has the same number of children $d$ (for $d \geq 2$), 
the toric ideal is generated in degree 2.
\end{conj}
\begin{conj}
There exists a quadratic Gr\"obner basis for the toric ideal of a binary tree.
\end{conj}

Using the Gr\"obner Walk \cite{CKM} implementation in {\tt magma},
we have computed thousands of Gr\"obner bases for random term orders
for the smallest binary trees.  It doesn't seem to be possible to
compute the entire Gr\"obner fan for these examples with {\tt CaTS} \cite{cats},
but the random computations have yielded some information: Conjecture 2 is tree for
the binary tree with 5 nodes, in fact, there are at least 4
distinct quadratic Gr\"obner bases for this tree.  Analysis of these
bases lends some optimism towards Conjecture 2.  However, for the
binary trees on 7 nodes, computation of over 1000 Gr\"obner bases did
not find a quadratic basis.  The best basis found contained
quartics and some bases even contained relations of degree 29.

Another nice family of toric ideals is given by $I_T$ for $T$ a path
of length $n$.  
Table~2 presents data for Markov bases of paths that 
leads us to conjecture that
this family also has well behaved ideals.
\begin{conj}
The toric ideal corresponding to a path is
generated in degree 3, with $2n - 4$ generators of
degree 3 needed.
\end{conj}

\begin{table}
\centering
\begin{tabular}{l|l|l|l|l}
\# of   & Degree   & \#Minimal  & Max  & Number \\
nodes & of $I_T$ & Generators   & degree  & of deg 3\\ \hline
3 & 6 & 6 & 3 & 2\\
4 & 19 & 32 & 3 & 4\\
5 & 36 & 102 & 3 & 6\\
6 & 61 & 259 & 3 & 8\\
7 & 90 & 540 & 3 & 10\\
8 & 127 & 1041 & 3 & 12\\
9 & 168 & 1842 & 3 & 14\\
10 & 217 & 3170 & 3 & 16\\
11 & 270  & 5286 & 3 & 18
\end{tabular}
\caption{Degree of $I_T$, size of Markov basis, maximum degree of a minimal generator, and
number of degree 3 generators for paths}
\label{tab:pathmarkov}
\end{table}



Unfortunately, the toric ideal of a general tree doesn't
seem to have such simple structure.  For $n \leq 9$, the
trees with highest degree minimal generators are those of the form
\psfig{file=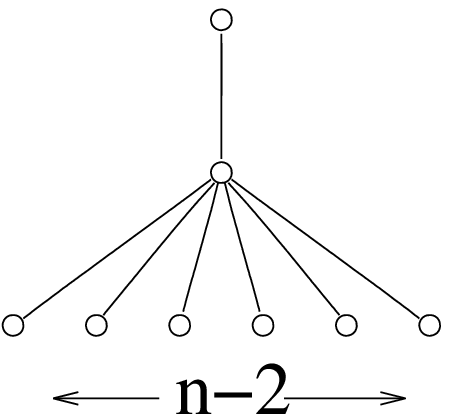,height=.4in}.  These trees 
require generators of degree $n$.

\section{Polytopes}\label{sec:poly}


In this section, we are interested in the following problem.
Given any observation  $(i_1, \dots , i_n)$ of the tree, 
which matrices $A = (a_{ij})$ make $p_{i_1 \dots i_n}$ 
maximal among the coordinates of the distribution $p$? 

To solve this problem, transform to logarithmic coordinates $b_{ij} =
\log(a_{ij})$.  Then the condition that $p_{i_1 \dots i_n} > p_{l_1
\dots l_n}$ for all $(l_1, \dots l_n) \in \{0,1\}^n$ is translated
into the the linear system of inequalities
\[ 
b_{i_1i_2} + \dots + b_{i_{\rho(n)}i_n} > 
    b_{l_1l_2}  + \dots + b_{l_{\rho(n)}l_n} 
\]
for all $(l_1, \dots l_n) \in \{0,1\}^n$.
The set of solutions to these inequalities 
is a polyhedral cone. For most values of $i_1,
\dots, i_n$, this cone will be empty.  Those sequences $i_1, \dots,
i_n$ for which the cone is maximal are called {\em Viterbi}
sequences. 
The collection of the cones, as
$(i_1, \dots, i_n)$ varies, is the normal fan of the polytope 
$P_T$, where $P_T$ is the convex hull of the columns of $\sA_T$.

Notice that $P_T$ is a polytope
in $\mathbb R^4$.  However, since all the monomials in (\ref{eq:toricdef}) 
are of degree $n-1$, we see that this polytope is actually contained
in $n-1$ times the unit simplex in $\mathbb R^4$.  Thus, $P_T$ is actually
a 3 dimensional polytope.
We call $P_T$ the {\em Viterbi} polytope.

The polytopes $P_T$ show remarkable finiteness properties as $T$ varies.
Since $P_T$ is defined as the convex hull of $2^n$ vectors, it would seem
that it could have arbitrarily bad structure.  However, as it is contained
in $n-1$ times the unit simplex, it can be shown that there are at most 
$O(n^{1.5})$ integral points in $P_T$.

\begin{example}\label{ex:path}
Eric Kuo has shown \cite{kuo} that if $T$ is a path with $n$ nodes, then $P_T$
has only two combinatorial types for $n > 3$, depending only on the parity of $n$.  
The polytope for the
path with 7 nodes is shown in 
Figure~1.
Think of this
picture as roughly a tetrahedron with the  vertex corresponding
to all $0 \to 1$ transitions and the vertex  
with all $1 \to 0$ transitions both sliced off (since if a path
has a $0 \to 1$ transition it must have a $1 \to x$ transition).
\end{example}

\begin{figure}
\centering
\epsfig{file=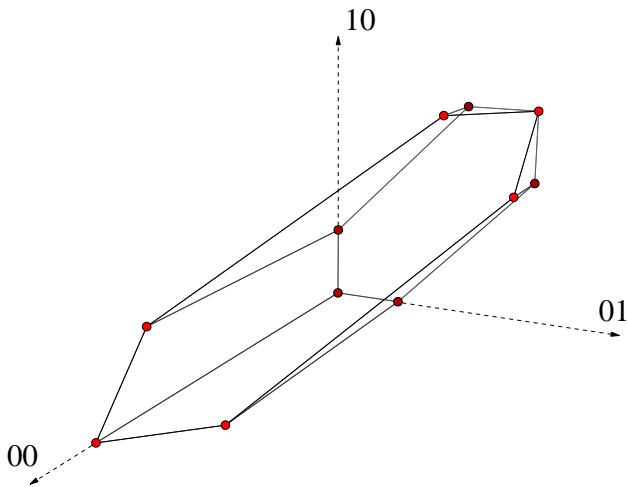, height=2.5in}
\caption{$P_T$ for $T$ a path with 7 nodes, after projecting
onto the first three coordinates $(b_{00},b_{01},b_{10})$.}
\label{fig:path7}
\end{figure}


Two facts from Example~\ref{ex:path} are important to remember.
First, the structure of the polytope is related more to the
topology of the tree than the size of the tree.  Second, there is a
distinction between even and odd length paths.  
We call a binary tree {\em completely odd} if the tree has all leaves
at an odd distance from the root. For example, 
the tree \psfig{file=binary4-2.ps,height=.35in} is completely odd.

\begin{theorem}\label{thm:1}
Let $T$ be a completely odd binary tree with more than three nodes.  
The associated polytope
$P_T$ always has the same combinatorial type with 8 vertices and 6
facets (see 
Figure~2).
\end{theorem}

\begin{proof}
First, we derive six inequalities that are satisfied by any binary
tree, deriving a ``universal'' polytope for binary trees.
Then we show that a completely odd binary tree has
labelings that give us all vertices of the ``universal'' polytope.

Thinking of the polytope space as the log space of the parameters $a_{ij}$, 
we write $\mathbb R^4$ with coordinates $b_{00}, b_{01}, b_{10}, b_{11}$.
Since $P_T$ lies in $n-1$ times the unit simplex in $\mathbb R^4$, we
have $b_{00} + b_{01} + b_{10} + b_{11} = n - 1$ and the 4
inequalities $b_{ij} \geq 0$.
We claim that any binary tree $T$ satisfies two additional inequalities
\begin{align}\label{eq:ineq}
\frac {b_{00} - b_{01}} {2} + b_{10} &\leq \frac {n + 1} 2,\\
\frac {b_{11} - b_{10}} {2} + b_{01} &\leq \frac {n + 1} 2.
\end{align}
We prove (\ref{eq:ineq}), the second inequality follows by
interchanging 1 and 0.

Fix a labeling of the binary tree.  We claim that the left hand side
of (\ref{eq:ineq}) counts the number of zeros that are ``created''
while moving down the tree, that is, it counts the number of leaves
that are zero minus one if the root is labeled zero.  Pick a non-leaf
of the tree which is labeled ``0''.  It has two children.  If both are
``0'', then this node contributes 2 to $b_{00} - b_{10}$.  If both are
``1'', then this node contributes -2 to $b_{00} - b_{10}$.  If one is
``0'' and one is ``1'', then the node doesn't contribute.  We think of
a ``0'' node with two ``0'' children as having created a new zero and
a ``0'' node with two ``1'' children as having deleted a zero.
Therefore we see that the term $(b_{00} - b_{10})/2$ counts the number
of zeros created as children of ``0'' nodes.  Similarly, if a
non-leaf is labeled ``1'', then its contribution to $b_{10}$ counts
the number of new zeros in the children.

Since there are $\frac {n+1} {2}$ leaves in a binary tree, there can
be at most $\frac {n+1} {2}$ zeros created, so (\ref{eq:ineq}) holds.
Notice that the labelings
that lie on this facet are exactly those with a one at the root and
all zeros at the leaves.


These six inequalities and the equality  
$b_{00} + b_{01} + b_{10} + b_{11} = n - 1$
define a three dimensional polytope in $\rr^4$.  We compute that there
are eight vertices of this polytope:
\begin{gather*}
(n-1, 0 ,0 ,0), \quad 
   (n-3, 0 ,2 ,0)\\
\left(\frac {n-3} 2 ,\frac {n+1} 2 ,0 ,0 \right ), \quad
   \left(0 ,\frac {2n} {3} ,\frac {n-3} {3} ,0 \right)\\
\left(0 ,\frac {n-3} {3} ,\frac {2n} {3} ,0 \right), \quad
  \left(0 ,0 ,\frac {n+1} 2 , \frac {n-3} 2 \right)\\
(0 ,2 ,0 ,n-3),\quad
  (0 ,0 ,0 ,n-1)
\end{gather*}
Six of these vertices occur in any binary tree:  a tree with all zeros
gives the $(n-1,0,0,0)$ vertex,   a tree with a one at the root and zeros
elsewhere gives $(n-3,0,2,0)$, and a tree with ones at the leaves and zeros
elsewhere gives $(\frac {n-3} 2, \frac{n+1} 2, 0,0)$.  Interchanging 1 and 0 gives 
three more vertices.  However, the remaining two vertices aren't obtained by
all binary trees.

The vertex $(0, \frac {n-3} 3, \frac {2n} 3, 0)$ lies on the facet defined by
(\ref{eq:ineq}), so we know it must have a one at the root, all zeros
at the leaves, and the labels must alternate going down the tree since there
are no zero to zero or one to one transitions.
This means that this vertex is representable by a labeled tree if and only if
the tree has all leaves at an odd depth from the root.  Notice that this
implies that $n$ must be divisible by 3 for the tree to be completely odd.
Finally, if $n > 3$ is odd and divisible by 3, then $n \geq 9$ and one checks
that the eight vertices are distinct.

See 
Figure~2
for a picture of the polytope and a Schlegel
diagram with descriptions of the labelings on the facets and at
the vertices.
\end{proof}

\begin{figure*}
\centering
\epsfig{file=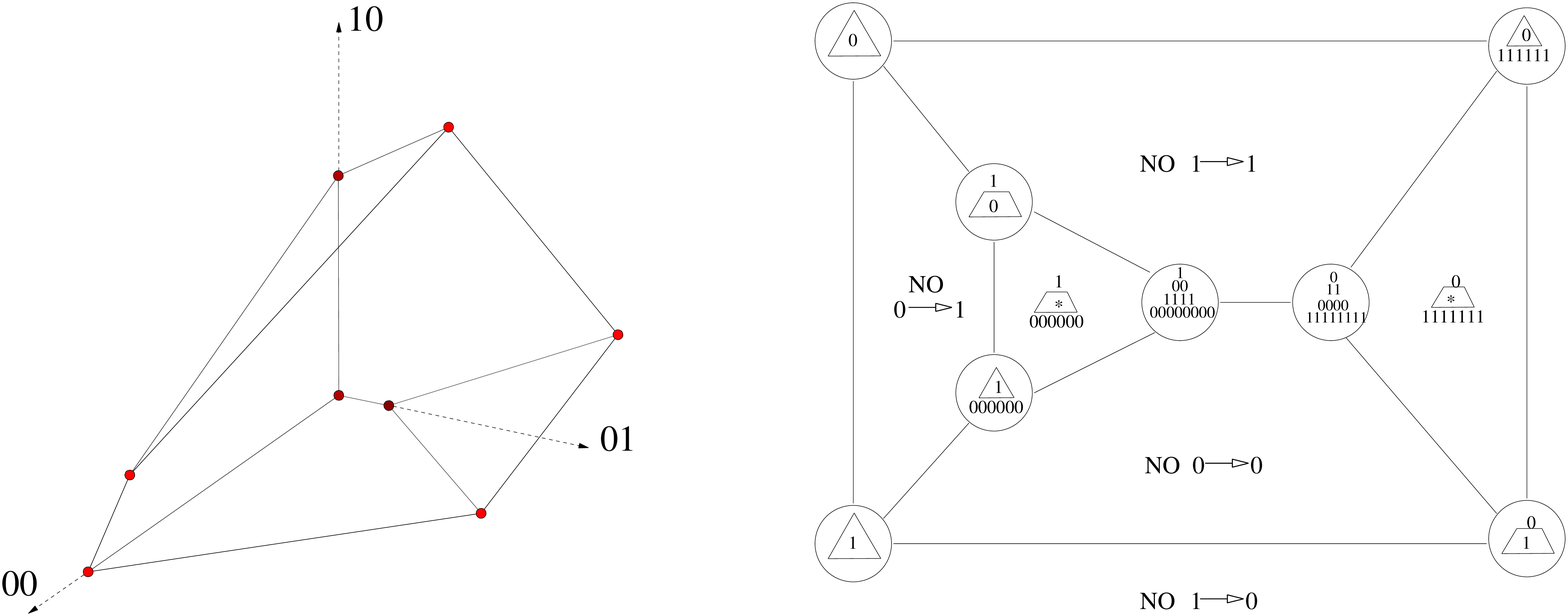, height=2.5in}
\caption{The polytope of the completely odd  
binary tree and a Schlegel diagram of this polytope
with facets and vertices labeled.}
\label{fig:schlegel}
\end{figure*}

In the case where $T$ is binary but not completely odd, the polytope shares 6
vertices with this universal polytope, but the remaining 2 vertices
are either not integral or not realizable.  However, the polytope
still shares much of the boundary with the universal polytope, so it is
perhaps realistic to expect that the polytope for a general binary
tree behaves well.  
Table~3
shows data from
computations for all binary trees with at most 23 nodes.  The maximum
number of vertices of $P_T$ appears to grows very slowly with the size
of the tree.
\begin{table}
\centering
\begin{tabular}{l|l|l|l|l}
Number &  Number of  & Min & Max & Ave\\
of nodes & binary trees & vertices & vertices & vertices \\\hline
3 & 1 & 4 & 4 & 4\\
5 & 1 & 7 & 7 & 7\\
7 & 2 & 8 & 10 & 9\\
9 & 3 & 8 & 13 & 11.33\\
11 & 6 & 10 & 14 & 11.66\\
13 & 11 & 11 & 13 & 11.91\\
15 & 23 & 8 & 16 & 14.35\\
17 & 46 & 12 & 17 & 13.82\\
19 & 98 & 10 & 20 & 14.65\\
21 & 207 & 8 & 19 & 14.8\\
23 & 451 & 10 & 20 & 15.6 
\end{tabular}
\caption{Minimum, maximum and average number of vertices of $P_T$ over all binary
trees with at most 23 nodes}
\label{tab:bintreepoly}
\end{table}

Although binary trees
seem to generally have polytopes
with few vertices, arbitrary trees are not so nice.  For example,
Figure~3 shows a tree with 15 nodes that has a
polytope with 34 vertices.

Table~4
shows data for all trees on at most 15 nodes.
It appears that the maximum number of vertices for the polytope
of an arbitrary tree of size $n$ grows approximately as $2n$.
Notice that the tree with all leaves at depth 1 has $P_T$ a tetrahedron,
giving the unique minimum number, 4, of vertices for all trees.

\begin{figure}
\centering
\epsfig{file=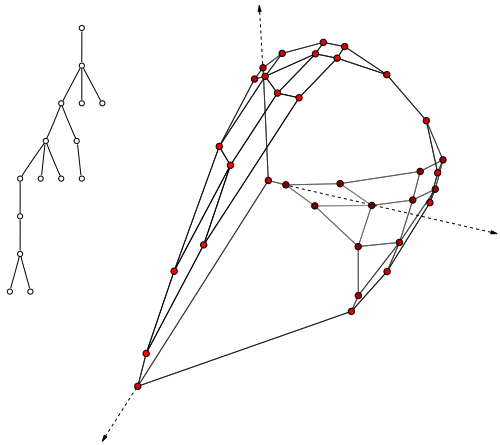, height=2.3in}
\caption{A tree $T$ with 15 nodes for which 
$P_T$ has 34 vertices, 58 edges, and
26 faces.}
\label{fig:badtreepoly}
\end{figure}

\begin{table}
\centering
\begin{tabular}{l|l|l|l|l}
Number &  Number of  & Min & Max & Ave\\
of nodes & trees & vertices & vertices & vertices \\\hline
3& 2& 4& 7& 5.5\\
4& 4& 4& 8& 7\\
5& 9& 4& 11& 8\\
6& 20& 4& 14& 9.7\\
7& 48& 4& 15& 10.75\\
8& 115& 4& 20& 12.59\\
9& 286& 4& 21& 13.67\\
10& 719& 4& 22& 15.42\\
11& 1842& 4& 25& 16.60\\
12& 4766& 4& 28& 18.3\\
13& 12486& 4& 31& 19.5\\
14& 32973& 4& 32& 19.75\\
15& 87811& 4& 34&  22.6
\end{tabular}
\caption{Minimum, maximum and average number of vertices of $P_T$ over all
trees with at most 15 nodes}
\label{tab:alltreepoly}
\end{table}

\begin{conj}
There is a bound on the number of vertices of $P_T$, where $T$ ranges
over all binary trees.  However, for an arbitrary tree, the number of vertices
of $P_T$ is unbounded.
\end{conj}

To extend these computations, a better algorithm for computing $P_T$
needs to be developed.  The naive algorithm for computing $P_T$
involves a loop of size $2^n$, elimination of duplicates points, and a
convex hull computation.  This algorithm can certainly be improved,
but it is not known whether there is a polynomial time algorithm for
constructing the polytope given a tree.  Is there a fast algorithm
that, given a tree $T$ and a point of $\rr^4$, outputs whether that
point arises from a labeling of $T$?  If so, then $P_T$ could be
constructed by testing the $O(n^{1.5})$ points inside $n-1$ times the unit simplex.


\section{Acknowledgments}
The author is grateful for useful conversations with Seth Sullivant, Lior Pachter, 
and Bernd Sturmfels and for much help from Raymond Hemmecke in fine tuning {\tt 4ti2}
for these computations.  The polytopal computations were done with {\tt polymake}.  
The author was supported by a National Defense Science and Engineering Graduate Fellowship.

\bibliographystyle{abbrv}

\end{document}